\documentclass[a4paper,12pt]{article}
\usepackage{amsfonts}
\usepackage{amssymb}
\usepackage{bbm}
\usepackage{amscd}
\usepackage{mathrsfs}
\usepackage{amsmath,amssymb, amsthm}
\usepackage{graphicx}
\usepackage{color}
\usepackage{authblk}
\usepackage{subcaption}
\usepackage{indentfirst,latexsym,amssymb}
\usepackage[bookmarks=false]{hyperref}
\usepackage{tikz}

\input amssym.def
\input amssym.tex

\newtheorem{lem}{Lemma}[section]%
\newtheorem{theorem}[lem]{Theorem}%
\newtheorem{exam}[lem]{Example}%
\newtheorem{prop}[lem]{Proposition}%
\newtheorem{rem}[lem]{Remark}%

 \def\s{\sigma}

 \def\O{\Omega}

    \def\oc{\overline C}
\def\ob{\overline B} 
  
 \def\op{\overline P} 
  
 \def\ox{\overline X}  \def\o1{\overline 1}

\def\olg{\overline g} \def\ola{\overline a} 
 \def\olz{\overline z}

\def\o{\overline}   \def\olb{\overline b}
\def\di{\bigm|} \def\lg{\langle} \def\rg{\rangle} \def\lls{\overline \s}

\def\Aut{\hbox{\rm Aut\,}}  \def\Syl{\hbox{\rm Syl}}

 \def\AGL{\hbox{\rm AGL}} \def\GL{\hbox{\rm GL}}  \def\P\GL{\hbox{\rm P\GL}}
  \def\FF{\mathbb{F}}
 \def\char{\hbox{\rm char}}

\def\C{\hbox{\rm C}}

\def\o{\hbox{\rm o}}

\def\char{ \, {\rm char}\,}

\def\calm{\mathcal{M}}

\def\ZZ{\mathbb{Z}}  

\def\nd{\mathrel{\bigm|\kern-.7em/}} 
 \def\f{\noindent}
\def\qed{\hfill $\Box$} \def\demo{\f {\bf Proof}\hskip10pt}

\begin{document}

\begin{center}
{\bf\large
Skew-Morphisms of Elementary Abelian $p$-Groups}{\footnotemark}
\end{center}

%\bigskip

\begin{center}
Shaofei Du,  Wenjuan Luo{\footnotemark},  Hao Yu\\
\medskip
{\it {\small
School of Mathematical Sciences, Capital Normal University,\\
Bejing 100048, People's Republic of China\\
dushf@mail.cnu.edu.cn,  wenjuan2202@163.com, 3485676673@qq.com
}}

\vskip 3mm
Junyang Zhang
\medskip
{\it {\small
\\School of Mathematical Sciences, Chongqing  Normal University,\\
Chongqing 401331, People's Republic of China\\jyzhang@cqnu.edu.cn
}}
\end{center}

\footnotetext{This work was supported by
the National Natural Science Foundation of China (12071312 and 11971248).}
\footnotetext{Corresponding author.}

\begin{abstract}
  A skew-morphism  of a finite group $G$ is a permutation $\s$ on $G$ fixing the
  identity element, and for which there exists an integer-valued function $\pi$ on $G$ such
  that $\s(xy)=\s(x)\s^{\pi(x)}(y)$ for all $x,y\in G$. It has been known that  for a given skew-morphism $\s $ of $G$, the product of the left regular representation of $G$ with $\lg \s \rg$ forms a permutation group on $G$, called a skew-product group of $G$. In this paper, the skew-product groups of finite elementary abelian $p$-groups are investigated. Some properties, characterizations and constructions about them are obtained.

 \end{abstract}

\section{Introduction}
\label{sec:intro}
Throughout the paper all groups  are finite. A \textit{skew-morphism} of a group $G$ is a permutation $\s$ on $G$ fixing the
identity element, and for which there exists an integer-valued function $\pi$ on $G$ such that
$\s(gh)=\s(g)\s^{\pi(g)}(h)$ for all $g,h\in G$.
In the special case that
$\pi(g)=1$ for all $g\in G$, the skew-morphism $\s$ is an automorphism of $G$. Thus a
skew-morphism is a generalization of a group automorphism.

\smallskip
The investigation of skew-morphisms is at least  related to the following two topics.

\smallskip
(1) {\it  Group factorizations}: Use
$L_{G}:=\{L_{g}\mid g\in G\}$ to denote the left regular representation of $G$.  Then both $\s$
and $L_g$ are permutations on $G$. For any $g,h\in G$,  we have
\begin{equation*}
(\s L_{g})(h)=\s(gh)=\s(g)\s^{\pi(g)}(h)=(L_{\s(g)}\s^{\pi(g)})(h),
\end{equation*}
and so $\s L_{g}=L_{\s(g)}\s^{\pi(g)}$.
Therefore,  $\langle \s\rangle L_{G}\subseteq L_{G}\langle \s\rangle $. Since $|\langle \s\rangle L_{G}|=|L_{G}\langle \s\rangle|$,  we have
$\langle \s\rangle L_{G}=L_{G}\langle \s\rangle $, which implies  that $X:=L_{G}\langle\s \rangle$ is a permutation group on $G$, called the \textit{skew-product} of $L_{G}$ by $\s$,  see
\cite{CJT2016, ZD2016}. Since $\langle \s \rangle  $  is a point stabilizer of the transitive
permutation group $X$, it is core-free in $X$, meaning that there is no nontrivial normal subgroup of $X$ contained in $\lg \s\rg $.

Conversely, suppose that $X$ is any group admitting a factorization $X=GY$
with $G\cap Y=1$ and  $Y=\langle y\rangle$ being cyclic and core-free in $X$. Then for any $g\in G$, there exist a unique $g'\in G$ and a unique nonnegative integer $i$ less than the order of $\langle Y\rangle$ such that $yg=g'y^{i}$. Define a permutation $\s$ on $G$ by $\s(g)=g'$, and an integer-valued function $\pi$ on $G$ by $\pi (g)=i$. Then one may check that $\s$ is a skew-morphism of $G$ with power function $\pi$.

\smallskip
(2) {\it Cayley maps:}
The concept of skew-morphism was first introduced as an algebraic tool to investigate
regular {\it Cayley maps}~\cite{JS2002}. A Cayley map $\calm =\mathrm{Cay}(G,S,P)$ is an embedding
of a (simple, undirected and connected) Cayley graph of a group $G$ with generating set
$S$ into an orientable closed surface such that, at each vertex $g$ of $\calm $, the local
orientation $R_g$ of the darts $(g,gx)$ incident with $g$ agrees with a prescribed
cyclic permutation $P$ of the generating set $S$, that is, $R_g(g,gx)=(g,gP(x))$ for all
$g\in G$ and $x\in S$. The automorphism group $\Aut(\calm)$ of a Cayley map $\calm$ contains a
vertex-regular subgroup induced by left multiplication of the elements of $G$. If the
cyclic stabilizer of a vertex is transitive on its adjacent vertices, then the automorphism
group $\Aut(\calm)$ is regular on the darts of $\calm$. In this case the map $\calm$ is called a
\textit{regular Cayley map}. It was shown by Jajcay and \v{S}ir\'{a}\v{n} that a Cayley map $\calm$ is
regular if and only if $P$ extends to a skew-morphism of $G$, see \cite[Theorem~1]{JS2002}. Thus the problem
of determining all regular Cayley maps of a group $G$ is equivalent to the problem of determining
all skew-morphisms of $G$ containing a generating orbit which is closed under taking inverses.

\smallskip
Now we are ready to recall the studying history of skew-morphisms of groups.
An interesting and important problem in this area is a determination of the skew-morphisms of a given family of groups. The problem seems challenging because even skew-morphisms of the cyclic groups have not yet been completely determined. As for cyclic groups,  skew-morphisms associated
with regular Cayley maps of cyclic groups can be extracted from the results obtained in~\cite{CT2014,Kwon2013}. The coset-preserving skew-morphisms of cyclic groups have been completely determined by
Bachrat\'y and Jajcay in~\cite{BJ2014, BJ2016}.
In~\cite{KN2011},
 Kov\'acs and Nedela proved that under certain numerical conditions a skew-morphism of a cyclic
 group may be decomposed into a direct product of skew-morphisms
of cyclic groups of prime power order, then in \cite{KN2017} they determined
the skew-morphisms of cyclic $p$-groups where $p$ is an odd prime, while  that of the case $p=2$ has been recently determined by Hu and the first author \cite{DH}.
Recently, Hu, Kwon and Zhang gave a classification of skew morphisms of cyclic groups which are square roots of automorphisms \cite{DH}.
 For partial results about skew-morphisms of  dihedral groups,
the reader is referred to~\cite{CJT2016,KMM2013,KK2016,Zhang2015,Zhang20152,ZD2016}. The corresponding regular Cayley maps of dihedral groups have been recently determined by Kov\'acs and  Kwon \cite{KK2021}.

\smallskip
A group is said to be \emph{monolithic} if it has a unique minimal normal subgroup which is nonabelian. The skew-morphisms of   monolithic groups were determined by Bachrat\'y, Conder and  Verret, see \cite{BCV}.
In particular, the skew-morphisms of  non-abelian simple group are determined in that paper. This result has been extended to that of  nonabelian characteristically simple groups, see ~\cite{CDL}. In this paper, we shall study
the skew-morphisms of  abelian characteristically simple groups (elementary abelian $p$-groups $\ZZ_p^n$).

\smallskip
Given  the \emph{skew product group} $L_G\langle \sigma\rangle$ of $L_G$ by $\sigma$,  for the purpose of this paper,  we may define  the \emph{skew product group} $X:=G\langle \sigma\rangle$ of $G$ by $\s$
as follows: every element of $X$ is uniquely written as $g\sigma^{i}$ where $g\in G$ and $i$ is a positive integer less than the order of $\sigma$; for each pair of elements $a\sigma^{i}, b\sigma^{j}\in X$, we have $(a\sigma^{i})(b\sigma^{j})=
a\sigma^{i}(b)\sigma^{\sum_{k=0}^{i-1}\pi(\sigma^{k}b)+j}$. It is straightforward to check by using the definition of the skew-morphism $\s$ that $X$ is indeed a group with operation defined above.  Sometimes, we just say $X$ a skew-product group of $G$ for short.

\smallskip
Now we are ready to state the main results of this paper. The  notations used in this paper are listed in the beginning of Section \ref{sec:pre}. By $G_X$ in Theorem~\ref{main1}, we denote  the core of $G$ in $X$.

\begin{theorem}\label{main1}
  Let $G\cong \ZZ_p^n$ where $p$ is a prime.   Let $X=G\lg \s \rg $ be  the  skew-product group of  $G$ by a skew-morphism $\s$, with order $|\s|=k\cdot p^m\ge 2$, where $m\ge 0$ and $p\nmid k$.
   Let $P=G\lg \s^k\rg $.
  Then  $P\lhd X$, either $G\lhd X$ or $G=G_X\times B$ where  $B\cong \ZZ_p$. More precisely,   one of the following holds:
   \begin{enumerate}
\item[\rm(1)] $m=0$ or $p=2$:  $X=G\rtimes \lg \s\rg $;
  \item[\rm(2)] $p\ne 2$, $m\ge 1$, $k=1$:  either $X=G\rtimes \lg \s\rg $ or  $X=(G_X\rtimes \lg \s \rg)\rtimes B$;
\item[\rm(3)]  $p\ne 2$, $m\ge 1$, $k\ge 2$:  $X=G\rtimes \lg \s\rg $; or

   \hskip 4cm $G\lhd P$,  $X=((G_X\times B)\rtimes \lg \s^k\rg )\rtimes \lg \s^{p^m}\rg $; or

   \hskip 4cm  $G\ntriangleleft P$, $X=((G_X\rtimes \lg \s^k\rg)\rtimes B)\rtimes \lg\s^{p^m}\rg $.
\end{enumerate}
\end{theorem}
In  Sections \ref{sec:zp2} and \ref{sec:zp3}, the number of skew-morphisms of $\ZZ_p^2$ and $\ZZ_p^3$ will be  determined,  respectively.
\begin{theorem} \label{main2} Let $X=G\lg \s\rg$ be a  skew-product group of $G\cong \ZZ_p^n$, where $n\in \{2, 3\}$. Then
\begin{enumerate}
  \item[\rm(1)]    $X\leq \AGL(\ell, p)$;
   \item[\rm(2)]  the number of skew-morphisms of $G$ is
 \begin{enumerate}
   \item[\rm(2.1)] $p=2$:\,  $1$, $6$ and $168$,  if $n=1, \,2,$ and $3$, respectively.
\item[\rm(2.2)] $p\ne 2$:\,  $p-1$, $2(p+1)(p-1)^3$ and $(p^3-1)(p^2-1)(p-1)(2p^3 -3p^2 +p+2)$,  if $n=1, \,2,$ and $3$, respectively.
\end{enumerate}
\end{enumerate}
\end{theorem}

\begin{rem}
 In  Theorem~\ref{main1}, some examples for  three  cases   when  $G\ntriangleleft X$ do exist.
 In Examples~\ref{e1}, ~\ref{e2} and ~\ref{e3},   minimal examples are respectively  constructed such that
 \vskip 3mm
   $X=(G_X\rtimes \lg \s \rg)\rtimes B$, where $G_X\cong \ZZ_3^4$, $|\s|=9$; or
 \vskip 3mm
    $G\lhd P$ and   $X=((G_X\times B)\rtimes \lg \s^k\rg )\rtimes \lg \s^{p^m}\rg $, where $G_X\cong \ZZ_3$ and $|\s|=6$; or
 \vskip 3mm
   $G\ntriangleleft P$ and  $X=((G_X\rtimes \lg \s^k\rg)\rtimes B)\rtimes \lg\s^{p^m}\rg $, where $G_X\cong \ZZ_3^4$ and $|\s|=18$.
\end{rem}

\begin{rem}
The statement for $G\lhd X$ when $m=0$ or $p=2$ was proved in   \cite{CJT2016} by using  some technical formulas  on  skew-morphisms.  In the present paper, we still cover it by using different methods from abstract group theories. Parts of results when $G=\ZZ_p^2$ in Theorem~\ref{main2}  were  also proved in  \cite{CJT2016}.
For the completeness  we reprove them by a different method too.
\end{rem}

After this introductory section, some  preliminary results will be given in Section \ref{sec:pre}. Theorem \ref{main1} will be proved in Section \ref{sec:proof}.
In Section \ref{sec:further}, some further properties  will be given and several examples related to Theorem~\ref{main1} are constructed.
In Sections \ref{sec:zp2} and \ref{sec:zp3},  the skew-morphisms of $\ZZ_p^2$ and $\ZZ_p^3$ will be determined and enumerated, respectively, that is Theorem~\ref{main2}.

\section{Preliminaries}
\label{sec:pre}
We at first list  some notations, to make the  paper more readable for researchers working on  Group theory, Combinatorics and other related areas.
\begin{enumerate}
\setlength{\itemsep}{0ex}
\setlength{\itemindent}{-0.5em}
  \item[] $\FF_q$ and $\FF_q^*$: the  field of order $q$ and its multiplicative group, resp.;
  \item[] $E_k$: identify matrix of degree $k$;
  \item[]  $|G|$ and  $|g|$:   the order of a group $G$ and  an element $g$ in $G$, resp.;
  \item[]  $H\leq G$ and $H<G$:  $H$ is a subgroup of $G$ and $H$ is a proper subgroup of $G$, resp.;
  \item[]  $H\lhd G$ and $H\char~G$:  $H$ is a normal and   characteristic subgroup of $G$, resp.;
  \item[] $G'$ and $Z(G)$: the derived subgroup and the center of $G$ resp.;
  \item[]  $G_X$: the core of $G$ in $X$ which is the maximal normal subgroup of $X$ contained in $G$;
  \item[] $M\rtimes N:$   a semidirect product of $M$ by $N$, in which $M$ is  normal;
  \item[] $C_K(H):=\{x\in K\di [x,H]=1\}$ and $N_K(H):=\{x\in K\di H^x=H\}$,  where $H, K\le G$;
  \item[]  $[a,b]:=a^{-1}b^{-1}ab$, the commutator of $a$ and $b$ in a group $G$;
  \item[]  $[a_{1},a_{2},\ldots,a_{n}]:
      =[[a_{1},a_{2},\ldots,a_{n-1}],a_{n}]$ where $n\geq 3$;
  \item[]  $[ia, jb]:=[a, b,  \overbrace{a, \cdots, a}^{i-1},   \overbrace{b, \cdots, b}^{j-1}]$ where $i, j$ are  positive integers;
  \item[] $\O_1(G)$: the subgroup $\lg g\in G\di g^p=1\rg$ of $G$ where $G$ is a $p$-group;
  \item[] $\Phi(G)$: the Frattini subgroup of $G$.
\end{enumerate}

\begin{rem}\label{derived} Recall that the Frattini subgroup $\Phi(G)$ of $G$ consists of elements $x$ of $G$ such that $G=\langle S,x\rangle$ leads to $G=\langle S\rangle$ for any $S\subseteq G$. It is well known that $G'\leq \Phi(G)$ if $G$ is a $p$-group. In particular, for a $p$-group $G$  and $H<X$, we have $X\neq HX'$.
\end{rem}
In what follows, we introduce some known group theoretical results for later uses.
A group $G$ is said to be {\it meta-abelian} if its derived group is abelian.

\begin{prop}[\cite{Xu1984}] \label{comm1}
Let $G$ be a meta-abelian group. Then, for $a,b\in G$ and a positive integer $n$, we have
$$(ab^{-1})^n=a^n(\prod_{i+j\leq n}[ia,jb]^{\binom{n}{i+j}})b^{-n}.$$
\end{prop}

\begin{prop}[\cite{Ito1}] \label{core0}
Let $G=AB$, where $A$ and $B$ are abelian subgroups of $G$. Then
\begin{enumerate}
  \item[\rm(1)]  $G$ is meta-abelian, that is, $G'$ is abelian;
  \item[\rm(2)]  if $G\ne 1$, then  $A$ or $B$ contains   a normal subgroup $N\ne 1$ of $G$.
\end{enumerate}
\end{prop}

\begin{prop}[\cite{Luc}] \label{cyclic}
If $G$  is a transitive permutation group of degree $n$  with a cyclic point-stabilizer,  then $|G|\le n(n-1)$.
 \end{prop}

The following result on meta cyclic $p$-groups  is easy to check directly.
\begin{prop} \label{metap}  Suppose that $G=\lg a\rg \rtimes \lg b\rg \cong \ZZ_{p^n}\rtimes \ZZ_p$, where $p$ is an odd prime. Then
$\O_1(G)=\lg a^{p^{n-1}}\rg \times \lg b\rg$. \end{prop}
\section{Proof of Theorem~\ref{main1}}
\label{sec:proof}
Throughout the remainder of the paper,
 we fix  the following notation in this section.
 Let $G\cong \ZZ_p^n$ and $C=\lg \s\rg \cong \ZZ_{kp^m}$, where  $(k,p)=1$.
Set $z=\s^{kp^{m-1}}$ if $m\ge 1$. Let $X=GC$ such that  $G\cap C=1$ and $C_X=1$. Let $P$ a Sylow $p$-subgroup of $X$ containing $G$.

\smallskip  Clearly, $P=P\cap (GC)=G(P\cap C)=G\lg \s^k\rg$  so that $X=PC_2$.  Then we first   need to deduce  the normality of $P$ in $G$.
\begin{lem}\,   \label{normal}  $P\lhd X$. \end{lem}
\demo   Since $X=GC$ where $G$ and $C$ are abelian groups, it follows from Proposition~\ref{core0} that $X$ is meta-abelian, that is $X'$ is abelian. Let $H$ be the Hall $p'$-subgroup of $X'$, which means that $|X':H|$ is a $p$-power. Then it follows from
$H\char X'\lhd X$ that $H\lhd X$. Since $HC_2\le X$ and $C_2$ is an Hall $p'$-subgroup of $X$,  we get
  $H\le C_2\le C$. Since $C_X=1$ and $H\lhd X$, we get $H=1$ and so $X'$ is a $p$-group. Since  $X'P$ is a $p$-subgroup where $P\in \Syl_p(X)$, we get   $X'\le P$, implying  $P\lhd X$.
\qed

\smallskip
In the following lemma,  we deal with a special: $k=1$, so that  $X=P$, a $p$-group.  For a $p$-group $P$, set  $\O_1(P)=\lg g\in P\di g^p=1\rg $. Remind that the  Frattini subgroup $\Phi(X)$ of $X$ consists of all non-generated elements of $X$.
\begin{lem}\label{p1} Suppose that $k=1$ and $m\ge 1$. Then  $\O_1(X):=G\rtimes \lg z\rg $ and one of the following holds:
\begin{enumerate}
  \item[\rm(1)] $p=2$:   $G\lhd X$;
  \item[\rm(2)] $p\ne 2$:  either $G\lhd X$ or $X=(G_X\rtimes \lg \s\rg )\rtimes B$, where $n\ge 4$  and $B\cong \ZZ_p$ is a complement of $G_X$ in $G$.
\end{enumerate}
\end{lem}
\demo  Let $M=GX'$. Then $M\lhd X$ and $M=M\cap X=G(M\cap C)$. Set $C_0=M\cap C=\lg\s_0\rg,$   where $|\s_0|=p^{m_0}$ for some $m_0\ge 1$.
By Remark~\ref{derived}, $M\lneqq X$.  Moreover, since  $(C_0)_M\lhd \lg G, C\rg=X$ and $C_X=1$, we get $(C_0)_M=1$.
In what follows we deal with two cases: $p=2$ and $p\ne 2$, separately.  Suppose that $z=\s^{p^{m-1}}$.

\smallskip
 {\it Case 1: $p=2$}

\smallskip
  Prove the lemma by the induction on $m$. Clearly, the result is true for $m=1$.  Assume $m>1$.
 By the inductive hypothesis, we get that $G\lhd M$ and so $M=G\rtimes C_0$.
  Let $H=\O_1(M)$.   Then  any involution $g$  in $M$ is of form $gz^i$ for some $g\in G$, which implies that
  $G\rtimes \lg z\rg\le H\le \lg G, z\rg$, that is $H=G\rtimes \lg z\rg$. It follows from  $H\char M\lhd X$ that $H\lhd X$. In particular,
   $H^\s=H$, that is $G^\s\le H=G\rtimes \lg z\rg$.

 \smallskip
   Suppose that  $G^\s \neq G$. Since $|H|=2^{n+1}$, we get  $H=GG^\s$ and
      $L:=G^\s \cap G\cong \ZZ_2^{n-1}$. Since $C_X=1$, we get $[z,G]\ne 1$, which means $H$ is nonabelian.  It follows from  $[L, G]=[L, G^\s]=1$ that $L=Z(H)$.
      Now we may write   $G^\s=L\times \lg gz\rg$ for some $g\in G\setminus L$, noting $H=GG^\s$.  Since $(gz)^2=g^2=z^2=1$, we get $[g, z]=1$, which implies  $[z, G]=[z, \lg L, g\rg]=1$,
      a contradiction.  Therefore, $G^\s=G$, that is $G\lhd X$.

\smallskip
 {\it Case 2: $p\ne 2$}
\smallskip
 Induction on $m$.  By the inductive hypothesis, we get that either $G\lhd M$ or $G_M\cong \ZZ_p^{n-1}$.
 Again, let $H=\Omega_1(M)$. Then $H\lhd X$.  First we prove a fact:

 \smallskip {\it Fact: $H=G\rtimes \lg z\rg $}

 \smallskip
{\it Proof of the fact:} If $G\lhd M$ then we  have  the completely same argument as  that in Case 1 to  see  $H=G\rtimes \lg z\rg$.
 Now suppose that $G_M\cong \ZZ_p^{n-1}$.    Set $G=G_M\times \lg b\rg$.
 Since  $M/G_M=\lg\overline\s_0\rg\rtimes \lg\overline b\rg \cong \ZZ_{p^{m_0}}\rtimes \ZZ_p$, where $m_0\ge 1$. By Proposition~\ref{metap},   we have
 $\O_1(M/G_M)=\lg \olz\rg \times \lg \olb\rg \cong \ZZ_p^2$. If $g\in M$ such that $|g|=p$, then $\olg\le \O_1(M/G_M)$ and  so $g\in \lg G_M, z, b\rg\le \lg G, z\rg $. Therefore,
  $\lg G, z\rg \le H\le \lg G,z\rg$, which implies $H=\lg G, z\rg $. From  $[\olz, \olb]=1$, we get $b^z\in bG_M$, that is $G^z=G$ and so $H=G\rtimes \lg z\rg $.

 \smallskip

 Come back to our proof.  Since $H=G\rtimes \lg z\rg $ and $H\lhd X$, we have either $G^\s=G$ or  $H=GG^\s$. In the second case,
   $Z(H)=G\cap G^\s\cong \ZZ_p^{n-1}$. It follows from  $Z(H)\char H\char M\lhd X$ that  $Z(H)\lhd X$. Therefore, $G_X=Z(H)\cong \ZZ_p^{n-1}$.

  \smallskip
   Finally, we show $|\O_1(X)=H.$  It is clearly true for   the case $G\lhd X$. Suppose $G\lneqq X$.  Then we have known that $X=(G_X\rtimes \lg \s\rg)\times \lg b\rg $ where $p\ne 2$.
   Now $X/G_X=\lg \overline \s\rg \rtimes \lg b\rg $. By Proposition~\ref{metap} again, every element of order $p$ in $X/G_X$ is of the form $(\overline \s)^i\olb^j$. Therefore,
    every element of order $p$ in $X$ is contained in $\lg G_X, b, z\rg=G\rtimes \lg z\rg =H$,that is $\O_1(X)\le H$, which in turn implies $\O_1(X)=H$.
 \qed

\vskip 3mm
Now we are ready to prove Theorem~\ref{main1}.

\vskip 3mm
\f {\bf Proof of Theorem~\ref{main1}:} From Lemma~\ref{normal}, $X=P\rtimes C_2$. In particular, if $m=0$, then $G=P\lhd X$, that is first part of Theorem~\ref{main1}.(1).
If $k=1$, then both a partial case  of (1) ($p=2$ and $k=1$) and  (2)   in Theorem~\ref{main1}  have been proved in Lemma~\ref{p1}. So we assume $k\ne 1$.

\smallskip Since a subgroup of $C_1$ is normal in $P=C_1G$ if and only if it is normal in $X=CG$, it follows from  $C_X=1$ that we get $(C_1)_P=1$, which implies
 $\s_1$ is a skew-morphism of $G$.    Therefore, $P=GC_1$ satisfies Lemma~\ref{p1} and in particular,  $H:=\O_1(P)=G\rtimes \lg z\rg\char P\lhd X$, which implies $H\lhd X$. So we deal with these two cases, separately.
\smallskip

\smallskip
 {\it (1) Case $p=2$:}

\smallskip
 Suppose that $G\ntriangleleft X$. Then $G\ne G^\s$ and $H=GG^\s$. As above, we have $Z(H)=G\cap G^\s$. By Lemma \ref{p1}, $P=G\rtimes C_1$. Note that $P=G\rtimes C_1=P^\s=G^\s\rtimes C_1$, hence $G^\s=Z(H)\times \lg gz\rg$ for some $g\in G\backslash Z(H)$. Since $(gz)^2=g^2=z^2=1$, we get $[g, z]=1$ and so $[z, G]=1$,  contradicting  to
 $C_X=1$. Therefore, $G^\s=G$, that is $G\lhd X$, which is second part of  Theorem \ref{main1}.(1).

\smallskip
{\it (1) Case $p\ne 2$:}

\smallskip
  Suppose that $G\ntriangleleft X$. Then $G\ne G^{\s}$ and by Lemma \ref{normal} we get $m\geq 1$. Again we have $H^{\s}=H$, $H=GG^\s$ and $Z(H)=G\cap G^{\s}\cong \ZZ_p^{n-1}$.
   Note that $Z(H)\lhd X$ as $Z(H)\char H\char P\char X$. Thus $G\ntrianglelefteq X$ but it contains a normal subgroup $Z(H)$ of $X$, so we get $G_X=Z(H)\cong \ZZ_p^{n-1}$. Set $G=G_X\times B$. If $G\lhd P$, then $X=((G_X\times B)\rtimes C_1)\rtimes C_2$;   if $G\ntriangleleft P$, then $X=((G_X\rtimes C_1)\rtimes B)\rtimes C_2$. Thus  Theorem \ref{main1}.(3) is proved.
   \vskip 3mm In summary, the proof of Theorem~\ref{main1} is completed.
\qed

\section{Further properties and constructions}
\label{sec:further}
An  more description for the case $P=(A\times C_1)\rtimes B$  is given in the following lemma,  where $B\cong \ZZ_p$.
\begin{lem}   \label{p^m1} Suppose that $P=(A\rtimes C_1)\rtimes B$, where $B=\lg b\rg\cong \ZZ_p$, $C_1\cong \ZZ_{p^m}$ and $m\ge 2$.  Viewing the conjugacy action of $\s$ on $A$ as a linear transformation,   we have $(\s-I)^{p^{m-1}-1}|_A\neq 0$, where $I$ is the identity.
 \end{lem}
 \demo First we give a  fact:

\smallskip
 {\it Fact 1:\, $C_P(A)=G\lg z\rg $, recalling $z=\s^{p^{m-1}}$ where $m\ge 2$}

\smallskip
 {\it Proof:}  Since $A\lneqq G$ and  $G\le C_P(A)\unlhd X$, we have $C_P(A)\ne G$. Then $G\lvertneqq C_P(A)=C_P(A)\cap (GC_1)=G(C_P(A)\cap C_1)$. Set $C_P(A)\cap C_1=\lg \s^{p^k}\rg $, which is non-trivial.
 In particular, $z\in C_X(A)$.  Since $\ox= X/A=\oc\ob=\oc \rtimes \ob$ where $B=\lg b\rg\cong \ZZ_p$,  we may set ${(\s^{p^{k}})}^b=\s^{lp^{k}}a$ for some  $a\in A$ and an integer $l$.
 Noting $\s^{p^k}\in C_X(A)$, we have   $$z^b=(\s^{p^{m-1}})^b=((\s^{p^k})^{p^{m-1-k}})^b=(\s^{lp^{k}}a)^{p^{m-1-k}}=z^la^{p^{m-1-k}}.$$
If $a^{p^{m-1-k}}=1$, then $z^b=z^l$ and $\lg z\rg \triangleleft X$, a contradiction. Therefore,  $a^{p^{m-1-k}}\ne 1$, which implies $m-1=k$ and so $C_P(A)\cap C_1=\lg z\rg.$ Therefore,
   $C_X(A)=G\lg z\rg$.

 \smallskip

 Return back to the proof of the lemma.  Since $|b|=p$, we may set $\s^b=\s^{p^{m-1}+1}a=z\s a$, where $a\in A$. By Proposition~\ref{core0},   $X$ is meta-abelian. It follows from Proposition~\ref{comm1}    that
 $$\begin{array}{lcl}
 z^b&=&(\s^{p^{m-1}})^b=(z\s a)^{p^{m-1}}=(\s a)^{p^{m-1}}\\
 &=&\s^{p^{m-1}}\prod_{i+j\le p^{m-1}} [i\s,\, ja^{-1}]^{\binom{p^{m-1}}{i+j}}a^{p^{m-1}} \\
  &=&z[(p^{m-1}-1)\s, a^{-1}].\end{array}$$

Now we show other fact:

\smallskip
{\it Fact 2:   $x(k):=[k\s, a^{-1}]=\prod_{i=0}^k(a^{(-1)^{k+i}\binom{k}{i}})^{\s^i}.$}

\smallskip
{\it Proof:} Now, $x(1)=[\s, a^{-1}]=a^\s a^{-1}=a^{-1}a^{\s}$.  By the induction hypothesis  on  $k$, we have
    $$ \begin{array}{lcl} &&x(k):=[k\s, a^{-1}]=[[(k-1)\s, a^{-1}], \s]\\
 &=&[(k-1)\s, a^{-1}]^{-1}[(k-1)\s, a^{-1}]^\s\\
 &=&\prod_{i=0}^{k-1}(a^{(-1)^{k+i+1}\binom{k-1}{i}})^{\s^i}\prod_{i=0}^{k-1}(a^{(-1)^{k+i}\binom{k-1}{i}})^{\s^{i+1}}\\
   &=&\prod_{i=0}^{k}(a^{(-1)^{k+i}\binom{k}{i}})^{\s^i}.$$
 \end{array}$$

Return back to the proof again.   View $\s$ as a linear transformation on $A$, that is $\s(a)=a^\s$, $a\in A$, a vector sapce. Since $\s$ is of order $p$-power, we know that that $\s-I$ is nilpotent.
  Then we have
 $$\begin{array}{lcl}
 z^b&=&zx(p^m-1)=z\prod_{i=0}^{p^{m-1}-1}(a^{(-1)^i(_{i}^{p^{m-1}-1})})^{\s^i}\\
     &=&z(\s-I)^{p^{m-1}-1}(a). \end{array}$$
If $x(p^m-1)=1$ as a group element (as a vector, $x(p^m-1)=0$), then $[z,b]=1$, forcing that $[z, G]=1$, a contradiction. So $x(p^m-1)\ne 1$, which gives $(\s-I)^{p^{m-1}-1}(a)\ne 0$ (as a vector). Therefore
  ${(\s-I)^{p^{m-1}-1}}|_A$ is not the zero transformation.
  \qed

\smallskip
Lemma~\ref{p^m1} helps us to construct some example when $G\lhd P$. See the following example.

\begin{exam} \label{e1} A  minimal example such that   $p$ is odd, $C=1$ and  $G\ntrianglelefteq P$.
 \end{exam}
{\it Construction:}  if $m=1$, then $G\lhd P$ and so take $m=2$.  By Lemma~\ref{p^m1},  $(\s-I)^{p-1}\ne 0$.  Since $\s-I$ is nilpotent, it is well-known that up to conjugacy, $\s$  can be taken a upper-triangle matrix.
 If $p$ is odd and $G=\ZZ_p^2$ or $\ZZ_p^3$, then $G\lhd P$.
So the minimal case when $p$ is odd  and $G\ntrianglelefteq X$  should be $p=3$ and $m=2$, so that $(\s-I)^2\ne 0$. Set
$$a_1=(1, 0, 0),\, a_2=(0, 1, 0), \, a_3=(0, 0, 1), M=\left(
                   \begin{array}{ccc}
                     1 & 1 & 0 \\
                     0 & 1 & 1 \\
                     0   &0  & 1 \\
                   \end{array}
                 \right),$$
 so that  $(a_1^\s, a_2^\s, a_3^\s)=(a_1, a_2, a_3)M$.   Turning back to the multiplication notation,
 $$a_1^\s=a_1, a_2^\s=a_1a_2, a_3^\s=a_2a_3.$$
 Let  $G=\lg a_1, a_2, a_3, b\rg\cong\ZZ_3^4$ and $C=\lg\s\rg\cong\ZZ_9$. Define a group $X$ as follows:
$$X=GC=\lg a_1,  a_2, a_3, b, \s\di a_1^\s=a_1, a_2^\s=a_1a_2, a_3^\s=a_2a_3,\s^b=\s^4a_3\rg .$$
Then we have $X=(A\rtimes C)\rtimes B$, where $A=\lg a_1, a_2, a_3\rg $  and $B=\lg b\rg$.

\smallskip
Suppose $X$ is an affine group. Then $X=T\rtimes H$, where $T\cong \ZZ_3^k$ and $H\le \GL(k, 3)$. Since $X$ contain an element of order 9 and both $\GL(2, 3)$ and
$\GL(3, 3)$ contains no elements of such order, we get $k \ge 4$.  It is easy to check that $X$ contains no normal subgroups of order $3^k$ for $k\ge 4$. Therefore, $X$ cannot be an affine group.
 \qed

\begin{exam} \label{e2} An minimal example when  $G\lhd P$ but $G\ntrianglelefteq X$.\end{exam}
{\it Construction:}\, Let  $$X=\lg a,b,\s|a^3=b^3=\s^6=1,a^{\s}=a^2b,b^{\s}=b^2\rg,$$
where $G=\lg a\s^2,b\rg$ and $P=\lg a, b, \s^2\rg $. Clearly,  $G\lhd P$ and   $G_X\cong \ZZ_p$.
Since $n=2$ and  $p=3$, this example is minimal one satisfying our condition.  \qed

\begin{exam} \label{e3}  An minimal  example when  $G\ntrianglelefteq P$.
 \end{exam}
{\it Construction:}\, Let
$$\begin{array}{lll} X&=&\lg a_1,a_2,a_3,a_4,a_5,\s|
\lg a_1,a_2,a_3,a_4,a_5\rg\cong\ZZ_3^5, \s^{18}=1,\\
&&
a_1^{\s}=a_1a_2,a_2^{\s}=a_2a_3,a_3^{\s}=a_3,a_4^{\s}=a_4,
\s^{a_5}=\s^{13}a_1a_2a_3\rg,\end{array} $$
where $G=\lg a_1,a_2,a_3,a_4,a_5\rg$ and $P=G\lg \s^2\rg $.
One can check that $G_P=\lg a_1, a_2, a_3, a_4\rg \cong \ZZ_3^4$ and $G\ntrianglelefteq P$ .
Clearly, this case cannot happen if $n\le 3$.  Checking by Magma, we know  that it cannot happen for $n=4$ too. So this example is minimal one.
\qed

\smallskip
 Finally, we give a sufficient condition for $X$  being an affine subgroup.
  \begin{lem}\label{AGL(n,p)}
 Suppose that $Z(P)\cong\ZZ_p^{n-1}$. Then  $X\le\AGL(n,p).$
 \end{lem}
\demo If $G\lhd X$, then $X=G\rtimes C\le\AGL(n,p),$ as desired. In what follows, we suppose $G\ntrianglelefteq X$.
Set $\s_1=\s^k$ and $\s_2=\s^{p^m}$ so that $C_1=\lg \s_1\rg $ and $C_2=\lg \s_2\rg .$
Set $P=GC_1$ and $Z=Z(P)$. Then  $\s_2$ normalizes $Z$. It is easy to see $\op=P/Z=\lg\ola,\lls_1\rg\cong \ZZ_{p^m}\rtimes\ZZ_p$~and $\lg\ola,\lls_1^{p^{m-1}}\rg~\char~\op$, where $a\in G\setminus Z$. Consider the conjugacy action of $C_2$ on $\lg\ola,\lls_1^{p^{m-1}}\rg\cong\ZZ_p^2$.  Since $(p, k)=1$ and $\s_2$ fixes $\lg \lls_1^{p^{m-1}} \rg$, we know that $\s_2$ fixes exactly two subgroups of order $p$ in $\lg\ola,\lls_1^{p^{m-1}}\rg$
(a  property of such elements in $\GL(2,p)$). So  there exists a $\lg \overline{a\s_1^{ip^{m-1}}}\rg $  fixed by  $\s_2$, for some integer $i$.   Then $\s_2$ fixes $T:=\lg Z, a\s_1^{ip^{m-1}}\rg\cong\ZZ_p^n $. Therefore,   $X=T\rtimes C\le \AGL(n,p) $.
 \qed

\section{Skew-morphisms of $\ZZ_p^2$}
\label{sec:zp2}
It is well-known that there are $p-1$ skew-morphisms of $\ZZ_p$, whose skew-product groups are $\ZZ_p\rtimes \ZZ_r$, where $r\di (p-1)$.
 In this section we shall determine skew-product groups of $\ZZ_p^2$ and the number of skew-morphisms of $\ZZ_p^2$.
 Suppose $p=2$. Then $X=\ZZ_2^2\rtimes \lg \s\rg$, where $\s\in \GL(2,2)\cong S_3$ and so  the number of skew-morphisms of $\ZZ_p^2$ is 6. In what follows, assume $p$ is odd.
  Copy Theorem~\ref{main2} when $G\cong \ZZ_p^2$ in here.
\vskip 3mm
\f {\bf Theorem~\ref{main2}}\ (Case $n=2$).
{ Let $X=G\lg \s\rg$ be a  skew-product group of $G\cong \ZZ_p^2$, where $p$ is odd. Then
  the number of  skew-morphisms of $G$ is  $2(p+1)(p-1)^3.$  Moreover,  $X\leq \AGL(2,p)$ and either $G\lhd X$ or $G_X\cong \ZZ_p$. }
\vskip 3mm
\demo    Let $X=GC$ where $G\cong \ZZ_p^2$,  $C=\lg \s\rg \cong \ZZ_{kp^m}$, $G\cap \lg \s\rg=1$ and $C_X=1$, where $(k,p)=1$.
Set $\s_1=\s^k$ and $\s_2=\s^{p^m}$ so that $\s=\s_1^u\s_2^v$ where $uk+vp^m=1$.
By Proposition~\ref{cyclic},  $m\in \{ 0, 1\}$.
  By Theorem~\ref{main1}, either $G\lhd X$ or $G_X\cong \ZZ_p$.  Let  $n$ be the number  of skew-morphisms of $G$. Then $n=n_1+n_2$, where
 $n_1$ and $n_2$  are  the numbers  of skew-morphisms of $G$ when $G\lhd X$ and $G \ntrianglelefteq X$, respectively.

\smallskip (1) Suppose that $G\lhd X$. Then  $\s$ may be any element of $\GL(2,p)$ and so  $n_1=p(p^2-1)(p-1)$.

\smallskip
(2) Suppose that $G\ntrianglelefteq X$. Then $m=1$.
Let $G\cong \ZZ_p^2.$ Then $P=G\rtimes \lg \s_1\rg  \le \AGL(2,p)$, where $\s_1$ is an element of order $p$.
 Let $M=\GL(2,p)$, where   $\lg\s_1\rg $ is  a Sylow $p$-subgroup of $M$. Then $C_M(\lg \s_1\rg )\cong \ZZ_p\rtimes (\ZZ_{p-1}\rtimes \ZZ_{p-1})$ and so $M$ contains the number $|M/C_M(\lg \s_1\rg)|=p+1$ of
 subgroups of order $p$. Therefore,     $\s_1$ has the number  $(p+1)(p-1)$ of choices.
In what follows, for each given $\s_1$,  we determine choices of $\s_2$ so that all choices of $\s=\s_1^u\s_2^v$ is known, where $uk+vp^m=1$.

\smallskip Let $P=G\rtimes \lg \s_1\rg $ for given $\s_1$ and $Z(P)=\lg z\rg $.  By Lemma~\ref{AGL(n,p)}, we get $X=T\rtimes \lg \s\rg$ where $T\cap G\cong Z(P)$ and moreover $\s_2$ normalizes an unique
 subgroup $T\cong \ZZ_p^2$ other than $\lg z,\s_1\rg $.  Then one may write   $T=\lg z, a\s_1^i\rg $ for some given $a\in G\setminus \lg z\rg $ and $i\in \ZZ_p^*$. Noting that different $i$ give different $T$ and so we have $p-1$
  such $T$. For each given $T$,  we shall determine the number of choices of $\s_2$ which normalizes it.
Identifying  $z$ with $(1,0)$ and $b\s_1^i$ with $(0,1)$, respectively so that $T=V(2,p)$,   the corresponding matrix representation  $P(\s_1)$ on $T$ is of the form $\left(\small\begin{array}{ccc}
    1  & 0 \\
    r  & 1 \\ \end{array}\right)$ for some given $r$, as $\s_1$ is of order $p$.  Since $[\s_2, \s_1]=1$,
    $P(\s_2)=sE_2$ for some $s\in \FF_p\setminus \{0, 1\}$.
Now,  the number of choices of $\s_1, i$ and $s$ are $p^2-1$, $p-2$ and $p-1$, respectively. Then we get
$n_2=(p^2-1)(p-2)(p-1)$, so that
     $$n=n_1+n_2=p(p^2-1)(p-1)+(p-2)(p+1)(p-1)^2=2(p+1)(p-1)^3.  $$  \qed

\smallskip  One of  examples with $G\ntrianglelefteq X$ is defined as follows:  let $P(\s_1)=\left(\small\begin{array}{ccc}
    1  & 1  \\
    0  & 1 \\ \end{array}\right)$, $P(\s_2)=sE_2$  where  $s\in \FF_p\setminus \{0, 1\}$, and   $i=1$.
Then $Z(P):=\lg t_{(1,0)}\rg $  and $G=\lg  t_{(1,0)}, t_{(0,1)}\s_1^{-1}\rg $.

\section{Skew-morphisms of $\ZZ_p^3$} \label{sec:zp3}
In this section we shall determine the skew-product groups  of $\ZZ_p^3$ and the number of skew-morphisms of $\ZZ_p^3$.

\smallskip
 Suppose $p=2$. Then $X=\ZZ_2^3\rtimes \lg \s\rg$, where $\s\in \GL(3,2)$ and so  the number of skew-morphisms of $\ZZ_p^3$ is 168. In what follows, assume $p$ is odd.
Before our discussion, we list some properties of $\GL(3,p)$ used later.
\begin{lem} \label{gl3p}
Let $L=\GL(3,p)$, where $p$ is an odd prime. Then $L$ has two conjugacy classes of elements  of order $p$, with the representatives $g_1$ and $g_2$,
 where $g_1$  fixes point-wise one 2-dimensional   subspace,  say $W_1$
   and $g_2$ cannot fix  point-wise any 2-dimensional   subspace but fix point-wise one  1-dimensional   subspace,  say $W_2.$
  Moreover,  $C_{L}(g_1)\cong \ZZ_p^3\rtimes \ZZ_{p-1}^2$ and it      fixes set-wise $W_1$.
   Let $\O$ be the set of  elements $x$ in $C_{L}(g_1)\setminus \{1\}$ such that $|x|\di (p-1)$ and $x$ moves every affine line $v+W$, where $v\ne 0$.
Then $|\O|=(p-2)(p-1)(p^2-p+1)$.
\end{lem}
\demo  Since the eigenvalues for every element $x$ of order $p$ is 1, there are two types of Jordan canonical matrices for $x$  with the properties stating in the  lemma, say $g_1$ and $g_2$, where
$$g_1=\left(\small\begin{array}{ccc}
                            1 & 1 &0 \\
                           0&1&0\\
                           0&0&1\\
                      \end{array}\right)  \quad {\rm and}\quad
                       g_2=\left(\small\begin{array}{ccc}
                            1 & 1 &0 \\
                           0&1&1\\
                           0&0&1\\
                      \end{array}\right).$$
  Direct checking shows that
  $$C_{L}(g_1)=\{\left(\small\begin{array}{ccc}
    x_1  & x_2 & x_3 \\
    0  & x_1 & 0 \\
    0 & x_4 & x_5\\
    \end{array}\right)\in L\di x_2,x_3,x_4 \in \FF_p,\, x_1,x_5\in \FF_p^*\}\cong \ZZ_p^3\rtimes \ZZ_{p-1}^2.$$
  Clearly,  $g_1$ fixes point-wise $W=\lg (0, 1, 0), (0, 0, 1)\rg $, which is fixed set-wise by $L$.
Let $g\in \O$. Since $g$ moves every affine line $v+W$, where $v\ne 0$, we have $x_1\ne 1$.
Since $|g|\di(p-1)$, one may check that either $g=x_1E_2$ where $x_1\in \FF_p\setminus\{ 0, 1\}$; or $g\in C_L(g_1)$ such that $x_1\neq x_5$, $x_2=\frac{x_3x_4}{x_5-x_1}$ and $x_3,x_4\in\FF_p$.
Thus,  $$|\O|=(p-2)+p^2(p-2)^2=(p-2)(p-1)(p^2-p+1).$$
\qed

Copy Theorem~\ref{main2} for $n=3$ in here.\vskip 3mm

\f {\bf Theorem~\ref{main2}:}  (Case $n=3$) {\it   Let $X=G\lg \s\rg$ be a  skew-product group of $G=\ZZ_p^3$, where $p$ is odd. Then
\begin{enumerate}
  \item[\rm(1)]    $X\leq \AGL(3,p)$; and either $G\lhd X$ or $G_X\cong \ZZ_p^2$;
  \item[\rm(2)]   the number of  skew-morphisms of $G$ is   $(p^3-1)(p^2-1)(p-1)(2p^3 -3p^2 +p+2)$.
\end{enumerate}}

\vskip 3mm The proof consists of  the following three lemmas.
\begin{lem}\label{m1} $m=1$. \end{lem}
\demo Let $X=GC$ where $G=\lg a,b,c\rg\cong \ZZ_p^3$, $C=\lg \s\rg \cong C_1\times C_2\cong \ZZ_{kp^m}$, where  $(k, p)=1$ and $m\in \{ 0, 1, 2\}$ (see Proposition~\ref{cyclic}).
  Set $\s_1=\s^k$, $\s_2=\s^{p^m}$, $\C_1=\lg \s_1\rg $ and $C_2=\lg \s_2\rg $.

\smallskip  Suppose that $G\lhd X$. Then  $X=G\rtimes \lg\s\rg\leq \AGL(3,p)$.

\smallskip Suppose that $G\ntrianglelefteq  X$.  Then by  Theorem~\ref{main1},  $m\ne 0$.
  Suppose that $m=2$.  Let $P\in \Syl_p(X)$.   By Theorem~\ref{main1},  $P\char  X$.
  Consider the group $P=GC_1$ which is order $p^5$, where $C_1$ is core-free. Suppose that $G\lhd P$. Since $C_1$ is core-free, we get
   $P\cong \ZZ_p^3\rtimes \ZZ_{p^2}\le \AGL(3, p)$, forcing that $\GL(3, p)$ contains   an element of order $p^2$, a contradiction.
  Therefore,   $G \ntrianglelefteq P$.
  Then by   Theorem~\ref{main1} we get $P=(A\rtimes C_1)\rtimes B$,
  where $A\cong \ZZ_p^2$, $C_1\cong \ZZ_{p^2}$ and $B\cong \ZZ_p$.
   From  $AC_1/C_{AC_1}(A)\lessapprox \Aut(A)\lessapprox \GL(2,p)$  where $AC_1$ is a $p$-group,
 we know that  $|C_{AC_1}(A):A|=1$ or $p$.
  Suppose that $|C_{AC_1}(A):A|=1$, that is  $[A, C_1]=1$. Then $C_1^p=(AC_1)^p\char AC_1\lhd P\char X$ and so $C_1^p\lhd X$, contradicting to $C_X=1$.   Suppose that $|C_{AC_1}(A):A|=p$. Then $[A, C_1]\ne 1$,  which implies $C_1^p=(C_{AC_1}(A))^p\char C_{AC_1}(A)\lhd GC_1$ and so $C_1^p\lhd \lg G, C\rg =X$, a contradiction too.  Therefore, $m=1$.
  \qed
\vskip 3mm

\begin{lem} \label{ag0}
  Suppose that $m=1$, $G\ntrianglelefteq X$ and $P\in \Syl_p(X)$.
  Then
  \begin{enumerate}
  \item[\rm(1)]    $P$ is a group of order $p^4$ defined by
$$P=\lg a, b, c\rg \rtimes \lg \s_1\rg \cong \ZZ_p^3\rtimes \ZZ_p,$$
where $\s_1$ has a faithful representation $P(\s_1)$ in $\GL(3,p)$ and it is just
 conjugate to the matrix $g_1$ in  Lemma~\ref{gl3p},  which fixes  pointwise a 2-dimensional subspace.
  \item[\rm(2)] $G_X\cong \ZZ_p^2$ and $X\le \AGL(3,p)$.
\end{enumerate}
\end{lem}
\demo
Since $m=1$, we have $X=P\rtimes \lg \s_2\rg$, $P=G\rtimes \lg \s_1\rg $ is of order $p^4$,   $k\le p^2-1$  and $C_1\cong \ZZ_p$.
 Identifying $G=\lg a, b, c\rg \cong \ZZ_p^3$  with 3-dimensional space,    $\s_1$  has  faithful representative $P(\s_1)$   in $\GL(3,p)$.
By  Lemma~\ref{gl3p},  $P(\s_1)$ fixes pointwise either  a 2-dimensional subspace or just a  1-dimensional subspace. If $P(\s_1)$ is in the second  case,  then  $G=\lg a, b, c\rg $ is the unique subgroup  isomorphic to $\ZZ_p^3$ of  $P$   and so it is   a characteristic subgroup of $P$, which implies $G\lhd X$, a contradiction.
 So $P(\s_1)$ must be in the first case, that is
 $$P=\lg a, b, c\rg \rtimes \lg \s_1\rg \cong \ZZ_p^3\rtimes \ZZ_p.$$
 Now $X=P\rtimes \lg\s_2\rg$, where  $\lg\s_2\rg$ acts faithfully on $P$ by conjugation. By the form of $P(\s_1)$, we know that  $Z(P)\cong \ZZ_p^2$.  It follows from  Lemma~\ref{AGL(n,p)}, $X=T\rtimes \lg \s\rg \le \AGL(3,p),$
\qed

\begin{lem}\label{number2}
There are $(p^3-1)(p^2-1)(p-1)(2p^3 -3p^2 +p+2)$ skew-morphisms of $\ZZ_p^3$.
\end{lem}
\demo
Let $n_1$, $n_2$ and $n$ have the same meaning as in  last section.

\smallskip
 Suppose that $G\lhd X$. Then  $\s$ can be taken any element of $\GL(3,p)$ and so  $n_1=p^3(p^3-1)(p^2-1)(p-1)$.

\smallskip
 Suppose that $G\ntrianglelefteq X$. Form the arguments of last lemma,
 $$P=\lg a, b, c\rg \rtimes \lg \s_1\rg \cong \ZZ_p^3\rtimes \ZZ_p,$$
where $\s_1$ has a faithful representation $P(\s_1)$ in $\GL(3,p)$,   which fixes point(vector)-wise the subspace, say $\lg b, c\rg$.
Then $P(\s_1)$ is conjugate to the matrix $g_1$ in  Lemma~\ref{gl3p}.
From the proof of Lemma~\ref{gl3p}, we know  $|C_{\GL(3,p)}(P(\s))|=p^3(p-1)^2$. Therefore,   we have the number $$|\GL(3,p)|/p^3(p-1)^2=(p^3-1)(p+1)$$
  of choices for $P(\s_1)$ and $\s$ as well.

\smallskip
 For each given $\s_1$, $Z(P)$ is uniquely determined and we set $Z(P)=\lg z_2,  z_3 \rg \le G$.  Similar to  the arguments in last section, $\s_2$ normalizes
 an unique group $T\cong \ZZ_p^3$ in $P$ other then $\lg Z(P), \s_1\rg$.
 So we write $T=\lg  z_1\s_1^i, z_2, z_3\rg $ for some $z_1\in G$ and $i\ne 0$. Then we determine the number of choices of $\s_2$ for each such $T$.
 Identifying $z_1\s_1^i,  z_2,$ and $z_3$ with $(1, 0, 0)$, $(0, 1, 0)$ and $(0,0,1),$ respectively and letting $T=V(3,p)$,
 we get another  representation $Q(\s)$ $\s$ in $\GL(3,p)$.  Under this base, we have
 $$Q(\s_1)=\left(\small\begin{array}{ccc}
    1 & r &s \\ 0&1&0\\ 0&0&1\\ \end{array}\right),\, r, s\in \FF_p, \, {\rm and }\, (r, s)\ne (0,0)$$
 and  $Q(\s_2)\in C_{\GL(3,p)}(Q(\s_1))$. Since $Q(\s_1)$ fixes  point-wise the subspace $W:=\lg (0, a_2, a_3)\di a_2, a_3\in \FF_p\rg $ and $[Q(\s_2), Q(\s_1)]=1$, we know that
 $Q(\s_2)$ fixes set-wise $W$,  and    does not fix $(1,0,0)$, equivalently,
 $Q(\s_2)$ moves every affine line $(k, 0, 0)+W$, where $k\ne 0$.

\smallskip
 By Lemma~\ref{gl3p}, we have   the number $(p-2)(p-1)(p^2-p+1)$  of choices for $\s_2$, relative to given $\s_1$ and $i$. By combining  the number of choices  of $\s_1$ and $i$, we get
   $$n_2=(p-2)(p-1)(p^2-p+1)(p^3-1)(p+1)(p-1)=(p^3-1)(p^2-1)(p^2-p-1)(p-1)(p-2);$$
$$\begin{array}{cll}
        n&=&n_1+n_2= p^3(p^3-1)(p^2-1)(p-1)+(p^3-1)(p^2-1)(p^2-p-1)(p-1)(p-2) \\
                &=&(p^3-1)(p^2-1)(p-1)(2p^3-3p^2+p+2).\\
        \end{array}$$ \qed
\vskip 3mm  One of  examples for $G\ntrianglelefteq X$ is the following:  let $\s_1=\left(\small\begin{array}{ccc}
    1 & 1 &0 \\
                           0&1&0\\
                           0&0&1\\ \end{array}\right)$, $\s_2=tE_3$, where $t\in\FF_p\setminus\{0,1\}$, and   $i=1$.
Then $Z(P)=\lg z_2,  z_3 \rg :=\lg t_{(0, 1, 0)},t_{(0,0,1)}\rg $  and $z_1=t_{(1, 0, 0)}\s_1^{-1}$ so that $G=\lg  t_{(0, 1, 0)},t_{(0,0,1)},t_{(1, 0, 0)}\s_1^{-1}\rg $.\qed

\end{document}